\documentclass[12pt]{amsart}
\usepackage{amssymb}
\usepackage{amsbsy}
\usepackage{amscd}
\advance\textheight by \topskip
\makeatletter
\renewcommand*\subjclass[2][2000]{%
  \def\@subjclass{#2}%
  \@ifundefined{subjclassname@#1}{%
    \ClassWarning{\@classname}{Unknown edition (#1) of Mathematics
      Subject Classification; using '2000'.}%
  }{%
    \@xp\let\@xp\subjclassname\csname subjclassname@#1\endcsname
  }%
}
\def\cal{\mathcal}
\def\Bbb{\mathbb}

\newenvironment{pf*}[1]{\proof[#1]}{\endproof}
\newcommand{\rom}{\textup}
\renewcommand{\thesubsection}{\thesection(\@roman\c@subsection)}
\makeatother
\newtheorem{Theorem}[equation]{Theorem}

\newtheorem{Lemma}[equation]{Lemma}
\newtheorem{Proposition}[equation]{Proposition}

\theoremstyle{definition}
\newtheorem{Definition}[equation]{Definition}
\newtheorem{Example}[equation]{Example}

\makeatletter
\renewcommand\section{\@startsection{section}{1}%
  {\z@}{.7\linespacing\@plus\linespacing}{.5\linespacing}%
  {\reset@font\normalfont\bfseries\centering}}
\makeatother

\theoremstyle{remark}
\newtheorem{Remark}[equation]{Remark}

\newtheorem*{Acknowledgements}{Acknowledgements}
\numberwithin{equation}{section}
\numberwithin{figure}{section}

\newcommand{\thmref}[1]{Theorem~\ref{#1}}
\newcommand{\secref}[1]{\S\ref{#1}}
\newcommand{\lemref}[1]{Lemma~\ref{#1}}
\newcommand{\propref}[1]{Proposition~\ref{#1}}

\newcommand{\defref}[1]{Definition~\ref{#1}}

\newcommand{\subsecref}[1]{\S\ref{#1}}

\newcommand{\exref}[1]{Example~\ref{#1}}

\newcommand{\Romnum}[1]{\expandafter\uppercase\expandafter{\romannumeral #1}} 
\newcommand{\C}{{\Bbb C}}
\newcommand{\Z}{{\Bbb Z}}

\newcommand{\R}{{\Bbb R}}
\newcommand{\CP}{\operatorname{\C P}}

\newcommand{\U}{\operatorname{\rm U}}
\newcommand{\SO}{\operatorname{\rm SO}}
\newcommand{\Spin}{\operatorname{\rm Spin}}
\newcommand{\Spinc}{\Spin^{c}}
\newcommand{\Map}{\operatorname{Map}}

\newcommand{\Ker}{\operatorname{Ker}}

\newcommand{\coker}{\operatorname{coker}}

\newcommand{\rank}{\operatorname{rank}}
\newcommand{\Sign}{\operatorname{Sign}}

\newcommand{\ind}{\mathop{\text{\rm ind}}\nolimits}

\newcommand{\M}{{\cal M}} % the moduli space
\newcommand{\G}{{\cal G}} % the gauge transformation group
\newcommand{\D}{{\cal D}} % 
\newcommand{\QQ}{{\cal Q}} % 
\newcommand{\LL}{{\cal L}} %
\newcommand{\E}{{\bar{E}}} %
\newcommand{\V}{{\bar{V}}} %
\newcommand{\f}{{\bar{f}}} %
\newcommand{\F}{{\cal F}} %
\newcommand{\VV}{{\cal V}} %
\newcommand{\W}{{\cal W}} %
\newcommand{\SW}{\operatorname{SW}}
\newcommand{\Zp}{\Z_p}
\begin{document}
\title[Mod $p$ vanishing theorem]{Mod $p$ vanishing theorem of Seiberg-Witten invariants for $4$-manifolds with $\Zp$-actions}
\author{Nobuhiro Nakamura}
\address{Research Institute for Mathematical Sciences, Kyoto university, Kyoto,
 	606-8502, Japan}
\email{nakamura@kurims.kyoto-u.ac.jp}
%\thanks{}
\begin{abstract}
We give an alternative proof of the mod $p$ vanishing theorem by F.~Fang of Seiberg-Witten invariants under a cyclic group action of prime order, and generalize it to the case when $b_1\geq 1$. 
Although we also use the finite dimensional approximation of the monopole map as well as Fang, our method is rather geometric. 
Furthermore, non-trivial examples of mod $p$ vanishing are given.
\end{abstract}

%57R57 Applications of global analysis to structures on manifolds, Donaldson and Seiberg-Witten invariants 
%57S17 Finite transformation groups
%
%
%\begin{abstract}
%\end{abstract}
%
\maketitle
%
%%%%%%%%%%%%%%%%%%%%%%%%%%%%%%%%%%%%%%%%%%%%%%%%%%%%%%%%%%%%%%%%%%%%%%%%%%%%%%%
\section{Introduction}\label{sec:intro}
%%%%%%%%%%%%%%%%%%%%%%%%%%%%%%%%%%%%%%%%%%%%%%%%%%%%%%%%%%%%%%%%%%%%%%%%%%%%%%%
%
In this paper, we investigate Seiberg-Witten invariants under a cyclic group action of prime order.
The Seiberg-Witten gauge theory with group actions has been studied by many authors \cite{Wang, Bryan, Fang0,Fang,RW,N, Cho, Bald, Sz} etc. 
Among these, we pay attention to a work by F.~Fang \cite{Fang}. 

In the paper \cite{Fang}, Fang proves that the Seiberg-Witten invariant of a smooth $4$-manifold $X$ of $b_1=0$ and $b_+\geq 2$ under an action of cyclic group $\Zp$ of prime order $p$, vanishes modulo $p$ if some inequality about the $\Zp$-index of Dirac operator and $b_+$ is satisfied, where $b_i$ is the $i$-th Betti number of $X$ and $b_+$ is the rank of a maximal positive definite subspace $H_+(X;\R)$ of $H_2 (X;\R)$. 
His strategy for proof is to use the finite dimensional approximation introduced by M.~Furuta \cite{Furuta} and appeal to equivariant $K$-theoretic devices such as the Adams $\psi$-operations.
This method requires concrete informations about equivariant $K$ groups.

On the other hand, in this paper, we give an alternative proof of Fang's theorem by a completely different method which is rather geometric. 
Then we are able to extend it to the case when $b_1 \geq 1$ by this geometric method.

To state the result, we need some preliminaries.

Let $G$ be the cyclic group of prime order $p$, and $X$ be a $G$-manifold. 
When $p=2$, we assume that the $G$-action is orientation-preserving. (Note that, when $p$ is odd, every $G$-action is orientation-preserving.) 
Fixing a $G$-invariant metric on $X$, we have a $G$-action on the frame bundle $P_{\SO}$. 
According to \cite{Fang}, we say that a $\Spinc$-structure $c$ is {\it $G$-equivariant} if the $G$-action on $P_{\SO}$ lifts to a $G$-action on the $\Spinc (4)$-bundle $P_{\Spinc}$ of $c$.

Suppose that a $G$-equivariant $\Spinc$-structure $c$ is given. 
Fix a $G$-invariant connection $A_0$ on the determinant line bundle $L$ of $c$. Then the Dirac operator $D_{A_0}$ associated to $A_0$ %together with the Levi-Civita connection of the fixed $G$-invariant metric becomes $G$-equivariant.
is $G$-equivariant, and 
the $G$-index of $D_{A_0}$ can be written as $\ind_{G} D_{A_0} = \sum_{j=0}^{p-1} k_j\C_j\in R(G)\cong\Z [t]/(t^p-1)$, where $\C_j$ is the complex $1$-dimensional weight $j$ representation of $G$ and $R(G)$ is the representation ring of $G$.

For any $G$-space $V$, let $V^G$ be the fixed point set of the $G$-action. Let $b_{\bullet}^G = \dim H_{\bullet}(X;\R)^G$, where $\bullet = 1,2,+$.
The Euler number of $X$ is denoted by $\chi(X)$, and the signature of $X$ by $\Sign (X)$.

In such a situation, F.~Fang \cite{Fang} proves the following theorem.
\begin{Theorem}[\cite{Fang}]\label{thm:Fang}
Let $G$ be the cyclic group of prime order $p$, and $X$ be a smooth closed oriented $4$-dimensional $G$-manifold with $b_1=0$ and $b_+\geq 2$. 
Let $c$ be a $G$-equivariant $\Spin^c$-structure. 
Suppose $G$ acts on $H_+(X;\R)$ trivially. 
If $2k_j\leq b_+ -1$ for $j=0,1,\ldots ,p-1$, then the Seiberg-Witten invariant $\SW_X(c)$ for $c$ satisfies
$$
\SW_X(c)\equiv 0 \mod p.
$$
\end{Theorem}

We will generalize \thmref{thm:Fang} to the case when $b_1\geq 1$. 
When $b_1\geq 1$, the whole theory can be viewed as a family on the Jacobian torus $J$. 
We consider the Jacobian torus $J$ as the set of equivalence classes of {\it framed} $\U(1)$-connections on $L$ whose curvatures are equal to that of the fixed $G$-invariant connection $A_0$. 
More concretely, $J$ is given as follows: 
Suppose that $X^G\neq \emptyset$, and choose a base point $x_0\in X^G$. 
Let $\G_0$ be the group of gauge transformations which are the identity at the base point $x_0$.  
Then the Jacobian $J$ is given as $J= (A_0 + i\ker d)/\G_0$, where $\ker d$ is the space of closed $1$-forms. 
Note that $G$ acts on $J$, and $J$ is isomorphic to $H^1(X;\R)/H^1(X;\Z)$ $G$-equivariantly. 

Since $J$ as above gives a well-defined family of connections, we can also consider the family of Dirac operators $\{D_A\}_{[A]\in J}$. 
Then its $G$-index $\ind_G \{D_A\}_{[A]\in J}$ is an element of the $G$-equivariant $K$-group $K_G(J)$ over $J$. 

Let $J^G = J_0\cup J_1\cup\cdots\cup J_K$ be the decomposition of the fixed point set $J^G$ into connected components. 
Choose a point $t_l$ in each $J_l$. For convenience, we assume that $J_0$ is the component including the origin which is represented by the fixed $G$-invariant connection $A_0$, and $t_0$ is the origin $[A_0]$.
By restriction, we have homomorphisms $r_l\colon K_G(J)\to K_G(t_l)$. 
Since each $K_G(t_l)$ is just the representation ring $R(G)\cong\Z[t]/(t^p-1)$, the image of $\alpha=\ind_G \{D_A\}_{[A]\in J}$ by $r_l$ is written as $r_l(\alpha) = \sum_{j=0}^{p-1}k_j^l\C_j$. 
(When $X^G=\emptyset$, a well-defined $G$-equivariant family of connections can not be constructed in general.
However coefficients $k_j^l$ can be defined ad hoc for our purpose.
See \subsecref{subsec:free}.) 
Now we state our main result which is a generalization of \thmref{thm:Fang}. 
%
%
%
%%%%%
\begin{Theorem}\label{thm:main}
Let $G$ be the cyclic group of prime order $p$, and $X$ be a smooth closed oriented $4$-dimensional $G$-manifold with $b_+\geq 2$ and $b_+^G\geq 1$. 
%Suppose that $X^G\neq \emptyset$ when $b_1\geq1$. 
Let $c$ be a $G$-equivariant $\Spinc$-structure, and $L$ be the determinant line bundle of $c$. 
Suppose $d(c) = \frac14 (c_1(L)^2 - \Sign (X)) - (1 - b_1 + b_+)$ is non-negative and even. 
If there exists a partition $(d_0,d_1,\ldots,d_{p-1})$ of $d(c)/2$ such that  $d_0 + d_1 + \cdots +d_{p-1} = d(c)/2$, and each $d_j$ is a non-negative integer and 
\begin{equation}\label{eq:ineq}
2k_j^l < 2d_j + 1 - b_1^G + b_+^G \quad (\text{for }j =0,1,\ldots ,p-1\text{ and any } l),
\end{equation}
then the Seiberg-Witten invariant $\SW_X (c)$ for $c$ satisfies 
\begin{equation*}\label{eq:mod-p}
\SW_X (c)\equiv 0 \mod p.
\end{equation*}  
\end{Theorem}
%%%%
%
%
%
\begin{Remark}
The number $d(c)$ is the virtual dimension of the Seiberg-Witten moduli space $\M_c$ of $c$, and $\SW_X (c)$ denotes the Seiberg-Witten invariant which is defined by the formula $\SW_X (c)= \langle U^{\frac{d(c)}2}, [\M_c]\rangle$, where $U$ is the cohomology class which comes from the $\U(1)$-action.
(See \defref{def:SW} below.)

When $b_1 > 0$, we can evaluate the fundamental class $[\M_c]$ by cohomology classes which originate in the Jacobian torus $J$ and define corresponding invariants.
Under our setting, there are some relations among these invariants which hold modulo $p$. 
This issue is treated separately in \secref{sec:torus}.  
\end{Remark}
\begin{Remark}
It can be easily seen that \thmref{thm:main} implies \thmref{thm:Fang}.
By the assumption of \thmref{thm:Fang}, $b_+^G=b_+\geq 2$ and $b_1=b_1^G=0$.
If $d(c)$ is odd or negative, then $\SW_X(c)=0$ by definition. 
Note that $d(c)$ is odd if and only if $b_+$ is even.
Therefore we can assume  $d(c)$ is non-negative and $b_+$ is odd. 
If the condition $2k_j\leq b_+-1$ for any $j$ is satisfied, then \eqref{eq:ineq} is satisfied for {\it any} partition of $d(c)/2$. 
Therefore we obtain \thmref{thm:Fang}. 
\end{Remark}
\begin{Remark}
\thmref{thm:main} can be rewritten in the following simpler form:
%The condition that there exists a partition $(d_0,d_1,\ldots,d_{p-1})$ of $d(c)/2$ such that  $d_0 + d_1 + \cdots +d_{p-1} = d(c)/2$, and each $d_j$ is a non-negative integer and satisfies \eqref{eq:ineq} can be rewritten as follows.
Let $X$ and $c$ be as in \thmref{thm:main}.  
Let $e_j$ (for $j=0,\ldots,p-1$) be integers defined by,
$$
e_j = \max_l\{(k_j^l -B),0\},
$$
where the constant $B$ is given as
$$
B=\left\{
\begin{aligned}
\frac12 &(1-b_1^G+b_+^G -1), \text{ when $1-b_1^G+b_+^G$ is odd},\\
\frac12 &(1-b_1^G+b_+^G -2), \text{ when $1-b_1^G+b_+^G$ is even}.
\end{aligned}
\right.
$$
%Then, the above condition is equivalent to $\sum_{j=0}^{p-1} e_j \leq d(c)/2$.
If $\sum_{j=0}^{p-1} e_j \leq d(c)/2$, then $\SW_X(c)\equiv 0\mod p$.
\end{Remark}
%%%%
Let us consider more precisely about lifts of the $G$-action to a $\Spinc$-structure. 
For a $\Spinc$-structure $c$, we have a bundle map $P_{\Spinc}\to P_{\SO}\times_X P_{\U(1)}$, where $P_{\U(1)}$ is the $\U(1)$ bundle for the determinant line bundle. 
This bundle map is a $2$-fold covering.
Suppose that $P_{\U(1)}$ is $G$-equivariant. 
If the action of a generator of $G$ on $P_{\SO}\times_X P_{\U(1)}$ lifts to $P_{\Spinc}$, then all of such lifts form an action on $P_{\Spinc}$ of an extension group $\hat{G}$ of $\Z_2$ by $G$:
\begin{equation}\label{eq:ext}
1\to\Z_2\to\hat{G}\to G\to 1.
\end{equation}

When $G$ is an odd order cyclic group, \eqref{eq:ext} splits. 
Therefore, if $\hat{G}$-lifts exists, then we can always take a $G$-lift on $P_{\Spinc}$. 
This is the case that $c$ is $G$-equivariant.

However, when $G=\Z_2$,  \eqref{eq:ext} does not necessarily split. 
The non-split case is when $\hat{G}=\Z_4$. 
In such a case, we say that the $\Z_2$-action is of {\it odd type} with respect to $c$. 
On the other hand, when $c$ is $\Z_2$-equivariant, we say that the $\Z_2$-action is of {\it even type} with respect to $c$.

Now suppose that the $\Z_2$-action is of {\it odd type} with respect to $c$.
For a $\Z_2$-connection $A$ on $L$, the Dirac operator $D_A$ is $\Z_4$-equivariant, and the $\Z_4$-index is of the form $\ind_{\Z_4} D_A = k_1 \C_1 + k_3 \C_3$. 
(This is because the $\Z_4$-lift of the generator of $\Z_2$ acts on spinors as multiplication by $\pm\sqrt{-1}$.)

In this case, we also have a result similar to \thmref{thm:main}. (Compare with Theorem 2 in \cite{Fang}.)
\begin{Theorem}\label{thm:odd}
Let $G =\Z_2$, and $X$ be a smooth closed oriented $4$-dimensional $G$-manifold with $b_+\geq 2$ and $b_+^G\geq 1$.
%Suppose that $X^G\neq \emptyset$ when $b_1\geq 1$, and that the $G$-action is of odd type with respect to a $\Spinc$-structure $c$.
Suppose that the $G$-action is of odd type with respect to a $\Spinc$-structure $c$.
For such $(X,c)$, \thmref{thm:main} holds as follows.
If there exists a partition $(d_1, d_3)$ of $d(c)/2$ such that  $d_1 + d_3 = d(c)/2$, and each $d_j$ is a non-negative integer and 
\begin{equation}\label{eq:ineq-2}
2k_j^l < 2d_j + 1 - b_1^G + b_+^G \quad (\text{for }j =1,3\text{ and any } l),
\end{equation}
then the Seiberg-Witten invariant $\SW_X (c)$ for $c$ satisfies 
\begin{equation*}\label{eq:mod-2}
\SW_X (c)\equiv 0 \mod 2.
\end{equation*}  
\end{Theorem}
Let us explain the outline of proofs of \thmref{thm:main} and \thmref{thm:odd}.

We also use a finite dimensional approximation $f$. %However, we treat $f$ more geometrically. 
We carry out the $G$-equivariant perturbation of $f$ to achieve the transversality, and then, under the assumption of \eqref{eq:ineq}, we see that the zero set of $f$ has no fixed point of the $G$-action by the dimensional reason concerning fixed point sets. 
Thus $G$ acts on the moduli space freely. 
Hence, if the dimension of moduli space is zero, then the number of elements in the moduli space is a multiple of $p$. 
From this, we can see that the Seiberg-Witten invariant is also a multiple of $p$. 
When the dimension of the moduli space is larger than $0$, it suffices to cut down the moduli space. 

To conclude the introduction, let us give a remark. 
At present, we did not find an application of \thmref{thm:main} in the case when $b_1\geq 1$. 
However, in the case of the $K3$ surface whose $b_1$ is $0$, the author and X.~Liu proved the existence of a locally linear action which can not be realized by a
smooth action by using the mod $p$ vanishing theorem \cite{LN}. 
Therefore, we could use \thmref{thm:main} or \thmref{thm:odd} to find such an action on a manifold with $b_1\geq 1$.
This problem is left to the future research.

The paper is organized as follows:
\secref{sec:finite} gives a brief review on the finite dimensional approximation of the monopole map and Seiberg-Witten invariants in the $G$-equivariant setting. 
\secref{sec:perturb} proves \thmref{thm:main} and \thmref{thm:odd}.
\secref{sec:torus} deals with Seiberg-Witten invariants obtained from tori in the Jacobian. 
\secref{sec:examples} gives some examples.
\begin{Acknowledgements}
The author would like to express his deep gratitude to M.~Furuta for invaluable discussions and continuous encouragements for years.
It is also a pleasure to thank Y.~Kametani for helpful discussions.
\end{Acknowledgements}
%%
%%%%%%%%%%%%%%%%%%%%%%%%%%%%%%%%%%%%%%%%%%%%%%%%%%%%%%%%%%%%%%%%%%%%%%%%%%%%%%%
\section{The $G$-equivariant finite dimensional approximation}\label{sec:finite}
%%%%%%%%%%%%%%%%%%%%%%%%%%%%%%%%%%%%%%%%%%%%%%%%%%%%%%%%%%%%%%%%%%%%%%%%%%%%%%%
The purpose of this section is to give a brief review on the finite dimensional approximation of the monopole map and Seiberg-Witten invariants in the $G$-equivariant setting.
%%%%
\subsection{The monopole map}
%%%%
Let $G=\Zp$, where $p$ is prime, and $X$ be a smooth closed oriented $4$-dimensional $G$-manifold with $b_+\geq 2$ and $b_+^G\geq 1$. 
Suppose that $X^G\neq \emptyset$. 
%(Later, we will see that it suffices to prove \thmref{thm:main} and \thmref{thm:odd} under this assumption, because the inequality \eqref{eq:ineq} can never be satisfied if $X^G= \emptyset$. See \lemref{lem:free}. )

Fix a $G$-invariant metric on $X$. 
Suppose a $\Spinc$-structure $c$ is $G$-equivariant.
We write $S^+$ and $S^-$ for the positive and negative spinor bundle of $c$.
Let $L$ be the determinant line bundle: $L=\det S^+$.

The Seiberg-Witten equations are a system of equations for a $\U(1)$-connection $A$ on $L$ and a positive spinor $\phi\in\Gamma (S^+)$,
%\begin{equation}\label{eq:SW}
%\left\{
%\begin{aligned}
% D_A \phi &= 0,\\
% F_A^+  &= q(\phi),
%\end{aligned}
%\right.
%\end{equation}
\begin{equation}\label{eq:SW}
\left\{
\begin{gathered}
 D_A \phi = 0,\\
 F_A^+  = q(\phi),
\end{gathered}
\right.
\end{equation}
where $D_A$ denotes the Dirac operator, $F_A^+$ denotes the self-dual part of the curvature $F_A$, and $q(\phi)$ is the trace free part of the endomorphism $\phi\otimes\phi^*$ of $S^+$ and this endomorphism is identified with an imaginary-valued self-dual $2$-form via the Clifford multiplication.

The action of the gauge transformation group $\G=\Map(X;\U(1))$ is given as follows:
for $u\in\G$, $u(A,\phi)=(A-2u^{-1}du, u\phi)$. 
Let $\M_c$ denotes the moduli space of solutions,
$$
\M_c=\{\text{solutions to \eqref{eq:SW}}\}/\G.
$$

Fix a $G$-invariant connection $A_0$ on $L$. 
Choose a base point $x_0$ in $X^G$, and let $\G_0 =\{u\in \G | u(x_0)=1\}$. 
Then $G$ acts on $\G_0$. 
The Jacobian torus $J$ is given as $J= (A_0+i\Ker d)/\G_0$, where $\Ker d$ is the space of closed $1$-forms. 

Let us define infinite dimensional bundles $\VV$ and $\W$ over $J$ by
\begin{align*}
\VV &= (A_0 + i\Ker d)\times_{\G_0}(\Gamma (S^+)\oplus \Omega^1(X)),\\
\W &= (A_0 + i\Ker d)\times_{\G_0}(\Gamma (S^-)\oplus \Omega^+(X)\oplus H^1(X;\R)\oplus\Omega^0(X)/\R),
\end{align*}
where $\R$ is the space of constant functions and $\G_0$-actions on spaces of forms and $H^1(X;\R)$ are trivial.
Note that $\VV$ decomposes into $\VV=\VV_\C\oplus\VV_\R$, where $\VV_\C$ is a complex bundle come from the component $\Gamma (S^+)$ on which $\U(1)$ acts by weight $1$, and $\VV_\R$ is a real bundle come from  $\Omega^1(X)$  on which $\U(1)$ acts trivially. 
The bundle $\W$  decomposes similarly as $\W=\W_\C\oplus\W_\R$.

To carry out appropriate analysis, we have to complete these spaces with suitable Sobolev norms. Fix an integer $k>4$, and take the fiberwise $L_k^2$-completion of $\VV$ and  the fiberwise $L_{k-1}^2$-completion of $\W$.  For simplicity, we use the same notation for completed spaces.

Now we define the monopole map $\Psi\colon\VV\to\W$ by
\begin{equation*}
\Psi (A,\phi, a) = (A, D_{A+ia}\phi, F_{A+ia}^+ - q(\phi), h(a), d^*a),
\end{equation*}
where $h(a)$ denotes the harmonic part of the $1$-form $a$. In our setting, $\Psi$ is a $\U (1)\times G$-equivariant bundle map.
Note that the moduli space $\M_c$ exactly coincides with $\Psi^{-1}(0)/\U (1)$.
%%%%
\subsection{Finite dimensional approximation}
%%%%
In this subsection, we describe the finite dimensional approximation of the monopole map according to \cite{Furuta-s}. 
(See also \cite{BF}.)

Decompose the monopole map $\Psi$ into the sum of linear part $\D$ and quadratic part $\QQ$, i.e., $\Psi = \D + \QQ$, where $\D\colon\VV\to\W$ is given by 
$$
\D(A,\phi,a)=(A,D_A\phi,d^+a,h(a),d^*a), 
$$
and $\QQ$ is the rest.

Let $W_\lambda$ (resp. $V_\lambda$) be the subspace of $\W$ (resp. $\VV$) spanned by eigenspaces of $\D\D^*$ (resp. $\D^*\D$) with eigenvalues less than or equal to $\lambda$. 
Let $p_\lambda\colon\W\to W_\lambda$ be the orthogonal projection.
As in \cite{Furuta}, we would like to consider $\D+p_\lambda \QQ$ as a finite dimensional approximation of $D + \QQ$. 
However  $W_\lambda$ and $p_\lambda$ do not vary continuously with respect to parameters in $J$.  
It is necessary to modify these.

Let $\beta\colon (-1,0)\to [0,\infty)$ be a compact-supported smooth non-negative cut-off function whose integral over $(-1,0)$ is $1$. 
For each $\lambda >1$, let us define the smoothing of the projection $\tilde{p}_\lambda\colon \W\to W_\lambda$ by
$$
\int_{-1}^0\beta(t)p_{\lambda+t}dt.
$$
Let $\iota_\lambda\colon W_\lambda\to\W$ be the inclusion.
Then the composition $\iota_\lambda\tilde{p}_\lambda$  varies continuously. 

For a fixed $\lambda$, we replace $W_\lambda$ with a vector bundle $W_f$ in the following lemma.
\begin{Lemma}[See \cite{Furuta-s}]\label{lem:Wf}
There is a $\U(1)\times G$-equivariant finite-rank vector bundle $W_f$ over $J$ and $\U(1)\times G$-equivariant bundle homomorphisms $\chi\colon W_f\to \W$ and $s\colon \W\to W_f$ which have the following properties.
\begin{enumerate}
\item The composition $\chi s$ on $W_\lambda$ is the identity. In particular, the image of $\chi$ contains $W_\lambda$. 
\item There is a $\U(1)\times G$-equivariant isomorphism from $W_f$ to the product bundle $J\times F_\C\oplus F_\R$, where $F_\C$ and $F_\R$ are complex and real representations of $G$ respectively. 
\end{enumerate}
\end{Lemma}
The proof of \lemref{lem:Wf} is given by modifying the proof of Lemma 3.2 in \cite{Furuta-s} $G$-equivariantly.

Let us consider the map $\D+\chi\colon \VV\oplus W_f\to \W$. 
Then we can show from \lemref{lem:Wf} that this map is always surjective. 
Therefore $V_f:= \Ker (\D+\chi)$ becomes a $\U(1)\times G$-equivariant finite-rank vector bundle. 

Now we can replace the family of linear maps $\D\colon V_\lambda\to W_\lambda$ with 
$$
\D_f\colon V_f\to W_f, \quad (v,e)\mapsto e,  
$$
which depends continuously on the parameter space $J$. 
Note that the formal difference $[V_f]-[W_f]$ gives the index of family $\D\colon V_\lambda\to W_\lambda$. 
In fact, it is easy to see that $\ker D \cong \ker D_f$ and $\coker D \cong \coker D_f$. 

For the non-linear part $\QQ$, we define a continuous family $\QQ_f\colon  V_f\to W_f$ by 
$$
\QQ_f (v,e) = -s\iota_\lambda\tilde{p}_\lambda\QQ (v).
$$

Then the map $\Psi_f :=\D_f +\QQ_f$ gives a finite dimensional approximation of $\Psi=\D +\QQ$ when we take sufficiently large $\lambda$. 
This is a $\U(1)\times G$-equivariant and proper map. In particular, the inverse image of zero is compact. 
%Furthermore $(V_f,W_f,\Psi_f)$ determines a stable homotopy element in some $\U(1)\times G$-equivariant stable homotopy group.
%
\begin{Remark}
The formulation in \cite{BF} is simpler than that of this section or \cite{Furuta-s}. 
However we need to use this formulation because the method in  \cite{BF} requires a trivialization of $\W$. 
In the non-equivariant setting, $\W$ can be always trivialized by Kuiper's theorem. 
However, in the $G$-equivariant setting, we do not know whether $\W$ can be trivialized $G$-equivariantly, or not.  
\end{Remark}
%
%%%%
\subsection{Seiberg-Witten invariants}
%%%%
%
%
Let $f_0=\Psi_f\colon V\to W$ be a finite dimensional approximation.
The space $V$ decomposes into the sum of a complex vector bundle $V_\C$ and a real vector bundle $V_\R$, $V=V_\C\oplus V_\R$, according to the splitting $\VV=\VV_\C\oplus\VV_\R$. 
Similarly $W =W_\C\oplus W_\R$.
Note that $[V_\C] - [W_\C]$ gives the $G$-index of the family of Dirac operators $\{D_A\}_{[A]\in J}$. 
Note also that $V_\R$ is a trivial bundle $\underline{F} = J\times F$, where $F$ is a real representation of $G$, and $W_\R = \underline{F}\oplus\underline{H}^+$, where $\underline{H}^+=J\times H^+(X;\R)$.

To obtain the Seiberg-Witten invariant, we need to perturb $f_0$ in general. 
For our purpose, we need to carry out the perturbation $G$-equivariantly. 
First, note that the moduli space $\M_c =f_0^{-1}(0)/\U(1)$ may have $\U(1)$-quotient singularities. 
(They are called {\it reducibles}. 
Strictly speaking, $f_0^{-1}(0)/\U(1)$ does not coincide with the genuine moduli space of solutions in general. 
However, after perturbation, the fundamental class of $f_0^{-1}(0)/\U(1)$ is equal to that of the perturbed moduli space.  
Therefore we abuse the term ``moduli space'' and the notation $\M_c$ for $f_0^{-1}(0)/\U(1)$.)
Let us consider the restriction of $f_0$ to the $\U(1)$-invariant part of $V$.
The $\U(1)$-invariant parts of $V$ and $W$ are $V^{\U(1)}=V_\R=\underline{F}$, and $W^{\U(1)}=W_\R = \underline{F}\oplus\underline{H}^+$, respectively.
%Note that the restriction $f_0|_{V^{\U(1)}}$ is a fiberwise linear proper map,thus, which  is just a fiberwise linear inclusion.
Since the restriction $f_0|_{V^{\U(1)}}$ is a fiberwise linear proper map, this is just a fiberwise linear inclusion.
Therefore, by fixing a non-zero vector $v\in H^+(X;\R)^G\setminus \{ 0\}$, and perturbing $f_0$ to $f = f_0 +v$, we can avoid reducibles, that is, $f^{-1}(0)^{\U(1)}=\emptyset$.
(Note that this perturbation is $\U(1)\times G$-equivariant.)

Let $\bar{V} = ((V_\C\setminus \{0\})\times_J V_\R)/\U(1)$, and define a vector bundle $\bar{E}\to\bar{V}$ by 
$$
\bar{E} = ((V_\C\setminus \{0\})\times_J V_\R\times_J W)/\U(1).
$$
Since $f$ is $\U(1)$-equivariant, $f$ induces a section $\bar{f}\colon\bar{V}\to\bar{E}$. 
Now, the moduli space $\M_c$ is the zero locus of $\bar{f}$.
Suppose $\bar{f}$ is transverse to the zero section of $\bar{E}$.
(In general, we need a second perturbation. Furthermore, in our case, the perturbation should be $G$-equivariant. This is a task in \secref{sec:perturb}.)
Then the moduli space $\M_c=\f^{-1}(0)$ becomes a compact manifold whose dimension $d(c)$ is
\begin{equation}\label{eq:dim-moduli} 
d(c) = \frac14 (c_1(L)^2 -\Sign(X)) - (1 -b_1 +b_+).
\end{equation}
We can determine the orientation of $\M_c$ from an orientation of $H^1(X;\R)\oplus H^+(X;\R)$.

Let us introduce a complex line bundle $\LL\to\bar{V}$ by $\LL=((V_\C\setminus \{0\})\times_J V_\R)\times_{\U(1)}\C$, where $\U(1)$ action on $\C$ is multiplication.
Let $U=c_1(\LL)$. 
Note that $H^*(\V;\Z)$ is isomorphic to $\Z [U]/(U^D-1)\otimes H^*(J;\Z)$ for some $D$ as an additive group. 

Now we give the definition of the Seiberg-Witten invariants. 
\begin{Definition}\label{def:SW}
The Seiberg-Witten invariant for a $\Spinc$-structure $c$ is given as a map,
$$
\SW_{X,c}\colon \Z [U] \otimes H^* (J;\Z)\to \Z,
$$
which is defined by $\SW_{X,c} (U^d\otimes\xi) = \langle U^d\cup \xi , [\M_c]\rangle$.
\end{Definition}
Note that an element $\xi$ in $H^* (J;\Z)$ can be written as a linear combination of Poincare duals of homology classes represented by subtori in $J$. 

Let $T$ be a subtorus in $J$, and its dimension be $d_T$. 
Suppose $d(c) - d_T$ is even and non-negative. 
Put $d^\prime = (d(c) - d_T)/2$. 
Then the Seiberg-Witten invariant $\SW_{X,c} (U^{d^\prime}\otimes P.D.[T])$ can be represented geometrically as follows: 
Let $\LL_1,\LL_2,\ldots ,\LL_{d^\prime}$ be $d^\prime$ copies of $\LL$ and $s_i\colon \V\to \LL_i$ $(i= 1,2,\ldots ,d^\prime)$ be arbitrary sections. 
Consider a section $\f_C$ of the vector bundle $\E\oplus\LL_1\oplus\cdots\oplus\LL_{d^\prime}$ given by $\f_C = (\f, s_1,\ldots ,s_{d^\prime})$.  
Now restrict $\f_C$ to $\V|_T$. 
If  $\f_C|_{\V|_T}$ is transverse to the zero section, then $\SW_{X,c} (U^{d^\prime}\otimes P.D.[T])$ is equal to the signed count of zeros of $\f_C|_{\V|_T}$ according to their orientations. 
(This method is called {\it cutting down} the moduli space.)

In this paper, we use the notation
$$
\SW_X(c) = \SW_{X,c}(U^{\frac{d(c)}2}),
$$
when $d(c)$ is non-negative and even.
%
%%%%%%%%%%%%%%%%%%%%%%%%%%%%%%%%%%%%%%%%%%%%%%%%%%%%%%%%%%%%%%%%%%%%%%%%%%%%%%%
\section{$G$-equivariant perturbation of $\bar{f}$}\label{sec:perturb}
%%%%%%%%%%%%%%%%%%%%%%%%%%%%%%%%%%%%%%%%%%%%%%%%%%%%%%%%%%%%%%%%%%%%%%%%%%%%%%%
%
In this section, we carry out the $G$-equivariant perturbation of $\f$, and finally prove \thmref{thm:main} and \thmref{thm:odd}.

Up to this point, we obtained a $G$-equivariant section $\f\colon\V\to\E$ which have no $\U(1)$-quotient singularity in the zero locus.
That is, the moduli space contains no reducible.
In order to go further, we need to identify $G$-fixed point sets $\V^G$ and $\E^G$. 
%
%%%%
\subsection{Fixed point sets $\V^G$ and $\E^G$}\label{subsec:fixed}
%%%%
%
Let us summarize the notation so far. 
The (perturbed) finite dimensional approximation is 
$$
f\colon V=V_\C\oplus \underline{F}\to W=W_\C\oplus \underline{F}\oplus\underline{H}^+.
$$ 
%where $\underline{F} = J\times F$, and $F$ is a real representation of $G$. 
The induced section is 
$$
\f\colon\V=(V_\C\setminus\{0\})/\U(1)\times_J \underline{F}\to\E=((V_\C\setminus\{0\})\times_J W_\C)/\U(1)\times_J (\underline{F}\oplus\underline{F}\oplus\underline{H}^+).
$$

Let us identify the fixed point set $\V^G=((V_\C\setminus\{0\})/\U(1)\times_J \underline{F})^G$. Note that $\V^G\to J^G$ is a fiber bundle. 
Recall that $[V_\C]-[W_\C]=\ind_G\{D_A\}_{[A]\in J}$. Then, for a fixed point $t_l\in J_l\subset J^G$, fibers of $V_\C$ and $W_\C$ over $t_l$ are written as
\begin{equation*}
V_\C|_{t_l} = \sum_{j=0}^{p-1}k_j^{l+}\C_j,\quad W_\C|_{t_l} = \sum_{j=0}^{p-1}k_j^{l-}\C_j,
\end{equation*}
and the relation $k_j^l = k_j^{l+} -k_j^{l-}$ holds.
Therefore the fiber of $\V^G$ over $t_l$ is $\V^G|_{t_l} = ((\sum_{j=0}^{p-1}k_j^{l+}\C_j\setminus\{0\})/\U(1))^G\times F_0$, where $F_0$ is the $G$-invariant part of the real representation $F$.
\begin{Lemma}\label{lem:proj}
There is a homeomorphism
$$ 
\left(\left(\left.\sum_{j=0}^{p-1}k_j^{l+}\C_j\setminus\{0\}\right)\right/\U(1)\right)^G\cong \coprod_{j=0}^{p-1} P(k_j^{l+}\C_j)\times \R_+,
$$
where $ P(k_j^{l+}\C_j)$ is the projective space of $k_j^{l+}\C_j$, and $\R_+$ is the set of positive real numbers.
\end{Lemma}
\proof
Note that there is a $G$-equivariant homeomorphism
$$
\left(\left.\sum_{j=0}^{p-1}k_j^{l+}\C_j\setminus\{0\}\right)\right/\U(1)\cong P(\sum_{j=0}^{p-1}k_j^{l+}\C_j)\times \R_+.
$$
A point $v$ in $P(\sum_{j=0}^{p-1}k_j^{l+}\C_j)$ is represented by a vector $(v_0,\ldots,v_{p-1})$ where $v_j\in k_j^{l+}\C_j$. 
Let $\zeta=\exp (2\pi\sqrt{-1}/p)$.
A point $v$ is fixed by the $G$-action if and only if there exists $\lambda\in\C\setminus\{0\}$ which satisfies $\lambda v_j = \zeta^j v_j$ for all $j$. 
Therefore there is a unique $j$ such that $v_j\neq 0$, and we have $\lambda =\zeta^j$ and $v_{j^\prime} = 0$ for all $j^{\prime}\neq j$.
Thus the lemma holds.
\endproof
By \lemref{lem:proj}, we see that $\V^G|_{t_l}\cong \coprod_{j=0}^{p-1} P(k_j^{l+}\C_j)\times \R_+\times F_0$. 
Therefore the dimension of the component $\V^G_{l,j}$ of $\V^G$ is given by
\begin{equation}\label{eq:dim}
\dim \V^G_{l,j} = 2k_j^{l+} -1 + a + b_1^G, 
\end{equation}
where $\V^G_{l,j}$ denotes the $j$-th component over $J_l\subset J^G$, and $a = \rank F_0$.
(Note that $b_1^G$ is the dimension of the base space $J_l$.)

Let us identify the fixed point set $\E^G$ similarly.
Note that 
\begin{align*}
\E&=((V_\C\setminus\{0\})\times_J W_\C)/\U(1)\times_J (\underline{F}\oplus\underline{F}\oplus\underline{H}^+) \\
\intertext{is an open submanifold of}
\E^\prime&:=((V_\C\oplus W_\C)\setminus\{0\})/\U(1)\times_J (\underline{F}\oplus\underline{F}\oplus\underline{H}^+).
\end{align*}
By the method similar to \lemref{lem:proj}, we see that $\E^{\prime G}|_{t_l}\cong\coprod_{j=0}^{p-1}P((k_j^{l+} +k_j^{l-})\C_j)\times\R_+\times (F_0\oplus F_0\oplus (H^+)^G)$.
Therefore the dimension of the component $\E^G_{l,j}$ of $\E^G$ is given by 
$$
\dim \E^G_{l,j} = 2(k_j^{l+} + k_j^{l-}) -1 + 2a + b_+^G + b_1^G, 
$$
where $\E^G_{l,j}$ denotes the $j$-th component over $J_l\subset J^G$.

Note that $\E^G\to\V^G$ is the disjoint union of vector bundles $\E_{l,j}^G\to \V_{l,j}^G$. 
The rank of $\E_{l,j}^G$ is given by
\begin{equation}\label{eq:rank}
\rank_\R \E_{l,j}^G =\dim \E^G_{l,j} - \dim \V^G_{l,j} = 2k_j^{l-} + a + b_+^G.
\end{equation}
%
%
%%%%
\subsection{Proof of \thmref{thm:main} in the case when $d(c)=0$.}\label{subsec:dim0}
%%%%
Suppose now that $d(c)=0$.
Under the assumption \eqref{eq:ineq}, formulae \eqref{eq:dim} and \eqref{eq:rank} imply that
$$
\dim  \V^G_{l,j} < \rank_\R \E_{l,j}^G.
$$
Therefore, we can perturb the section $\f\colon\V\to\E$ on a small neighborhood of the fixed point set $\V^G$ $G$-equivariantly so that $\f$ has no zero on $\V^G$.
Then it is easy to carry out a $G$-equivariant perturbation outside the $G$-fixed point sets so that $\f$ is transverse to the zero section. 
(For instance, consider on quotient spaces $\V/G$ and $\E/G$, and then pull back to original spaces.)

Note that the moduli space $\M_c =\f^{-1}(0)$ no longer contains any $G$-fixed point. 
Hence $G$ acts freely on $\M_c$. 
Thus we have $\SW_X(c)\equiv 0\mod p$.
%
%%%%
\subsection{Proof of \thmref{thm:main} in the case when $d(c)$ is positive and even}\label{subsec:dim-even}
%%%%
%
Let us introduce $G$-equivariant complex line bundles $\LL_j$ over $\V$ $(j=0,\ldots,p-1)$ by 
$$
\LL_j = ((V_\C\setminus\{0\})\times_J V_\R)\times_{\U(1)}\C_j,
$$
 and fix $G$-equivariant sections $s_j\colon\V\to\LL_j$. (It is easy to make a $G$-equivariant section. Choose an arbitrary non-$G$-equivariant section, and average it by the $G$-action.)
We will cut down the moduli space by these $(\LL_j,s_j)$.
 
%Suppose $d(c)$ is even. %, and let $\bar{d}=\frac12 d(c)$.
Fix a partition $(d_0,d_1,\ldots,d_{p-1})$ of $d(c)/2$ such that $d_j\geq 0$ and $d_0+d_1+\cdots +d_{p-1}=d(c)/2$. 
Instead of the section $\f\colon\V\to\E$, we consider 
\begin{align*}
\f_C&\colon\V\to\E\oplus d_0\LL_0\oplus\cdots\oplus d_{p-1}\LL_{p-1} =:\E_C\\
\intertext{which is defined by}
\f_C &= (\f,s_0,\ldots,s_0,s_1,\ldots,s_{p-1}).
\end{align*}
Hereafter, we argue in analogous way to that of \subsecref{subsec:fixed}. 
We write $(\E_C)^G_{l,j}$ for the component of the fixed point set $(\E_C)^G$ over $\V^G_{l,j}$. 
Then the rank of the vector bundle $(\E_C)^G_{l,j}\to\V_{l,j}$ is given by
\begin{equation}
\rank_\R (\E_C)^G_{l,j} = 2(k_j^{l-} + d_j)  + a + b_+^G. 
\end{equation}
An argument similar to that of the case when $d(c)=0$ in \subsecref{subsec:dim0} completes the proof of \thmref{thm:main} when $X^G\neq\emptyset$.

\subsection{The case when $X^G=\emptyset$}\label{subsec:free}
%To conclude this subsection, we give a remark about the case when $X^G=\emptyset$. 
The base point $x_0 \in X^G$ is used for the well-defined $G$-equivariant family of connections over the Jacobian $J$. 
When $b_1=0$, we do not need the base point to construct a finite dimensional approximation. 
Therefore, the argument in this section also works in the case when $b_1=0$ and $X^G=\emptyset$. 
%On the other hand, even though $b_1\geq 1$ and $X^G=\emptyset$, if we have a well-defined $G$-equivariant family of connections over $J$, then the $G$-index of the family of Dirac operators is defined.
On the other hand, in the case when  $b_1\geq 1$ and $X^G=\emptyset$, we can define coefficients $k_j^l$ ad hoc for our purpose, although we do not have a well-defined $G$-equivariant family of connections.
Consider the Jacobian $J$ as $J=(A_0+i\ker d)/\G$, where $\G$ is the {\it full} gauge transformation group. 
Decompose the $G$-fixed point set $J^G$ into connected components: $J^G=J_0\cup\cdots\cup J_K$. 
Choose a point $t_l$ in each component $J_l$ and a connection $A_l$ in each class $t_l$.   
We assume that $J_0$ is the component of $[A_0]$ and $t_0=[A_0]$, where $A_0$ is the fixed $G$-equivariant connection. 
Then, for each $A_l$, we can redefine the $G$-action on the $\Spinc$-structure $c$ such that $A_l$ is fixed by the redefined $G$-action. 
(This is proved as in \lemref{lem:action}.)
Then the Dirac operator $D_{A_l}$ is $G$-equivariant, and the $G$-index $\ind_G D_{A_l}$ is written as $\ind_G D_{A_l}=\sum_{j=0}^{p-1}k_j^l\C_j$.
In such a situation, we can prove the following.
\begin{Lemma}\label{lem:free}
Suppose that $d(c)$ in \eqref{eq:dim-moduli} is nonnegative and even. 
If $X^G=\emptyset$, then there is no partition $(d_0,d_1,\ldots,d_{p-1})$ of $d(c)/2$ which satisfies \eqref{eq:ineq}. 
\end{Lemma}
\proof
Coefficients $k^l_j$ are calculated by the $G$-index theorem. 
(See \subsecref{subsec:kjl}.)
In fact, we can show that
$$
k^l_0 = k^l_1= \cdots = k^l_{p-1} = \frac1p \ind D_{A_0} = \frac1{8p} ( c_1(L)^2-\Sign (X)),
$$
for any $l$.
Note that $1-b_1+b_+=p(1-b_1^G+b_+^G)$ when $X^G=\emptyset$. 
(This follows from the formulae $\chi (X) =p\chi(X/G)$ and $\Sign(X)=p\Sign(X/G)$.)  
Therefore \eqref{eq:ineq} is equivalent to $\frac1p d(c) < 2d_j$ for $j=0,1,\ldots,p-1$. 
Summing up these equations from $j=0$ to $p-1$ implies a contradiction. 
\endproof
Therefore, the assumption $X^G\neq \emptyset$ can be omitted logically.
%
%
%
%%%%
\subsection{Proof of \thmref{thm:odd}}
%%%%
%
Let $G=K=\Z_2$ and $\hat{G}=\Z_4$, and consider the short exact sequence,
$$
0\to K \to\hat{G}\to G\to 0.
$$
If the $G$-action is of odd type with respect to a $\Spinc$-structure $c$, then $\hat{G}$ acts on the whole theory.
In this case also, as in \secref{sec:finite}, we obtain the $\U(1)\times \hat{G}$-equivariant finite dimensional approximation
$$
f\colon V=V_\C\oplus \underline{F}\to W=W_\C\oplus \underline{F}\oplus\underline{H}^+.
$$ 
Note that the $\hat{G}$-action on $J$, $\underline{F}$ and $\underline{H}^+$ factors through the surjection $\hat{G}\to G$, and hence the actions of the subgroup $K\subset \hat{G}$ on $J$, $\underline{F}$ and $\underline{H}^+$ are trivial.

We need to identify $K$-fixed point sets as well as $\hat{G}$-fixed point sets. 
Note that $K$-actions on $V_\C$ and $W_\C$ are given as multiplication by $-1$ on each fiber, which are absorbed by $\U (1)$-actions. 
Therefore $K$-actions on $\V$ and $\E$ are trivial. 

Thus we see that the $\hat{G}$-action on the section $\f\colon\V\to\E$ is reduced to an action of $G=\hat{G}/K$. 
Then, an argument analogous to \subsecref{subsec:fixed},\subsecref{subsec:dim0}, \subsecref{subsec:dim-even} and \subsecref{subsec:free} proves \thmref{thm:odd}.
%
%
%
%
%%%%%%%%%%%%%%%%%%%%%%%%%%%%%%%%%%%%%%%%%%%%%%%%%%%%%%%%%%%%%%%%%%%%%%%%%%%%%%%
\section{Cutting down the moduli by tori in $J$}\label{sec:torus}
%%%%%%%%%%%%%%%%%%%%%%%%%%%%%%%%%%%%%%%%%%%%%%%%%%%%%%%%%%%%%%%%%%%%%%%%%%%%%%%
%
This section deals with Seiberg-Witten invariants obtained from tori in $J$.
In this section, let $G=\Z_p$ where $p$ is prime, and suppose that $X$ is a closed oriented $4$-dimensional $G$-manifold with $b_+\geq 2$, $b_+^G\geq 1$ and $b_1\geq 1$, and $X^G\neq\emptyset$. Let $c$ be a $G$-equivariant $\Spinc$-structure.

First, we suppose that a subtorus $T$ in $J$ is $G$-invariant, i.e., $T=gT$ for $g\in G$.  
Let $d_T=\dim T$.
Suppose that $d(c) - d_T$ is non-negative and even, and put $d^{\prime} = \frac12 (d(c) - d_T)$.
For a partition $(d_0,d_1,\ldots,d_{p-1})$ of $d^\prime$, consider $\f_C\colon \V\to \E_C= \E\oplus d_0\LL_0\oplus\cdots\oplus d_{p-1}\LL_{p-1}$ as in \subsecref{subsec:dim-even}.
Then consider the restriction $\f_C|_{\V|_T}$ of $\f_C$ to $\V|_T$.    
By perturbing $\f_C|_{\V|_T}$ $G$-equivariantly in the way similar to that of \secref{sec:perturb}, we can prove the following. 
\begin{Theorem}
Let $d_T^G=\dim T^G$. Suppose that $X^G\neq\emptyset$ and that $d(c) - d_T$ is non-negative and even. 
Put $d^{\prime} = \frac12 (d(c) - d_T)$.
If there exist a partition $(d_0,d_1,\ldots,d_{p-1})$ of $d^\prime$ such that  $d_0 + d_1 + \cdots +d_{p-1} = d^\prime$, and each $d_j$ is a non-negative integer and 
\begin{equation*}
2k_j^l < 2d_j + 1 - d_T^G + b_+^G \quad (\text{for }j=0,1,\ldots,p-1\text{ and any }l),
\end{equation*}
then 
\begin{equation*}
\SW_{X,c} (U^{d^\prime}\otimes P.D.[T])\equiv 0 \mod p,
\end{equation*}  
where $k_j^l$ are defined similarly from $\ind_G\{D_A\}_{[A]\in T}\in K_G(T)$.
\end{Theorem}

On the other hand, when $T$ is not $G$-invariant, the following holds.
\begin{Theorem}
Let $d_T^G=\dim T^G$. Suppose that $X^G\neq\emptyset$  and that $d(c) - d_T$ is non-negative and even. 
Put $d^{\prime} = \frac12 (d(c) - d_T)$.
If there exist a partition $(d_0,d_1,\ldots,d_{p-1})$ of $d^\prime$ such that  $d_0 + d_1 + \cdots +d_{p-1} = d^\prime$, and each $d_j$ is a non-negative integer and 
\begin{equation*}
2k_j^l < 2d_j + 1 - d_T^G + b_+^G \quad (\text{for }j=1,2,\ldots,p-1\text{ and any }l),
\end{equation*}
then 
\begin{equation*}
\sum_{i=0}^{p-1}\SW_{X,c} (U^{d^\prime}\otimes P.D.[g^iT])\equiv 0 \mod p.
\end{equation*}  
\end{Theorem}
\proof
Let us consider $\tilde{T} = T\cup gT\cup g^2T\cup \cdots\cup g^{p-1}T$ for $g\in G$, and  the restriction $\f_C|_{\V|_{\tilde{T}}}$ of $f_C$ to $\V|_{\tilde{T}}$.
Note that $\tilde{T}$ is not necessarily a manifold. 
Let $T_k$ be the set of $t\in\tilde{T}$ such that the number of $g^iT$ ($i=0,1,\ldots,p-1$) which contains $t$ is lager than or equal to $k$, that is,    
$$
T_k = \{t\in \tilde{T}\,|\, \#\{i\,|\,t\in g^iT\}\geq k\}.
$$
Note that $T_1 = \tilde{T}$ and $T_p = \bigcap_{i=0}^{p-1} g^i T$.
Then $\tilde{T}=T_1\supset T_2\supset \cdots\supset T_p$ gives a stratification. 
Note that $\dim\tilde{T}=\dim T_1 > \dim T_2$.
Note also that $T_p$ is $G$-invariant and contains all fixed points. 
By perturbing $\f_C|_{\V|_{T_p}}$ $G$-equivariantly in the way similar to \secref{sec:perturb}, $\f_C|_{\V|_{T_p}}$ comes to have no zero. 
(This is due to a dimensional reason.) 
Next perturb $\f_C$ on ${\V|_{T_{p-1}\setminus T_p}}$ $G$-equivariantly so that $\f_C|_{\V|_{T_{p-1}\setminus T_p}}$ has no zero. 
Successively perturb $\f_C$ on ${\V|_{T_{k}\setminus T_{k+1}}}$ for $k>1$ $G$-equivariantly so that $\f_C|_{\V|_{T_{k}\setminus T_{k+1}}}$ has no zero. 
Finally, carry out a $G$-equivariant perturbation of $\f_C|_{\V|_{\tilde{T}}}$ outside $V_{T_2}$ to achieve the transversality with the zero-section. 
Since all zeros are on ${\V|_{\tilde{T}\setminus T_2}}$, and $G$ acts freely on the set of zeros, the conclusion holds.
\endproof
%%
%
%
%
%%%%%%%%%%%%%%%%%%%%%%%%%%%%%%%%%%%%%%%%%%%%%%%%%%%%%%%%%%%%%%%%%%%%%%%%%%%%%%%
\section{Examples}\label{sec:examples}
%%%%%%%%%%%%%%%%%%%%%%%%%%%%%%%%%%%%%%%%%%%%%%%%%%%%%%%%%%%%%%%%%%%%%%%%%%%%%%%
%
The purpose of this section is to give several examples. 
In order to apply \thmref{thm:main} and \thmref{thm:odd} to concrete examples, we need to calculate coefficients $k_j^l$. 
Therefore we first discuss how to calculate coefficients $k_j^l$. 
%
%
%%%%%%
\subsection{How to calculate $k_j^l$}\label{subsec:kjl}
%%%%%%
%
%
Recall that we decomposed the fixed point set $J^G$ of the Jacobian torus into connected components: $J^G=J_0\cup\cdots\cup J_K$, and chose a point $t_l$ in each $J_l$. 
Fix a generator $g\in G$, and write $\hat{g}$ for the action of $g$ on the $\Spinc$-structure $c$. 
For the origin $t_0=[A_0]$, by definition, it holds that $\hat{g}A_0 = A_0$. 
Therefore, we can calculate $k_j^0$ by the $G$-index formula such as $\ind_g D_{A_0} = \text{(contributions from fixed points)}$. 
First we briefly review the $G$-index formula. 
(See \cite{AS2,AS3,AH,AB}.)

Let $X^G = X_0\cup X_1\cup\cdots\cup X_N$ be the decomposition of the fixed point set $X^G$ into connected components, where $X_0$ is assumed to be the component of the base point $x_0$. 
Then, the $G$-index formula for $t_0=[A_0]\in J^G$ is written as 
\begin{equation*}%\label{eq:g-ind}
\ind_g D_{A_0} =\sum_{j=0}^{p-1} \zeta^j k_j^0 = \sum_{n=0}^N \F^0_n(g), 
\end{equation*}
where $\zeta=\exp (2\pi\sqrt{-1}/p)$ and each $\F_n^0(g)$ is a complex number associated to the component $X_n$ which is given as follows.

Let $L_n$ be the restriction of the determinant line bundle $L$ to $X_n$.
Then $g$ acts on each fiber of $L_n$ as the multiplication with a complex number $\nu_n$ of absolute value $1$.
(In our case, $\nu_n$ is a $p$-th root of $1$.)

There are two cases with respect to the dimension of $X_n$. 
Since we assume the $G$-action is orientation-preserving, the dimensions of $X_n$ are even.

If $X_n$ is just a point $x_n$, the tangent space over $x_n$ is written as 
$$
T_{x_n}X = N(\omega_1)\oplus N(\omega_2),
$$
where $N(\omega_j)$ is the complex $1$-dimensional representation on which $g$ acts by multiplication with $\omega_j$. 
(In our case, $\omega_j$ is a $p$-th root of $1$.)

Then the number $\F_n^0(g)$ is given by,
\begin{equation}\label{eq:fn0}
\F_n^0(g) = \nu_n^{\frac12}\frac1{\omega_1^{1/2} - \omega_1^{-1/2}}\frac1{\omega_2^{1/2} - \omega_2^{-1/2}}.
\end{equation}
The right hand side is only defined up to sign. 
To determine the sign precisely, we need to see the $g$-action on the $\Spinc$-structure $c$. 
When $G$ is the cyclic group of odd order $p$ and the $\Spinc$-structure $c$ is $G$-equivariant, signs of $\omega_i^{1/2}$ and $\nu_n^{1/2}$ are determined by the rule that
\begin{equation}\label{eq:sign}
\left(\omega_i^{1/2}\right)^p = \left(\nu_n^{1/2}\right)^p  = 1.
\end{equation}
(See \cite[p.20]{AH}.)
On the other hand, when $p=2$, it is somewhat subtle problem to determine the sign precisely.
(See \cite{AB}.)

If $X_n$ is a $2$-dimensional surface $\Sigma_n$, the restriction of the tangent bundle of $X$ to $\Sigma_n$ is written as 
$$
TX|_{\Sigma_n} = T\Sigma_n\oplus N(\omega),
$$
where $N(\omega)$ is the normal bundle of $\Sigma_n$ in $X$, and $g$ acts on the fiber of $N(\omega)$ as multiplication with $\omega$.

In this case,  $\F_n^0(g)$ is given as,
\begin{equation}\label{eq:fn0-2}
\F_n^0(g) = -\nu_n^{\frac12}\cdot\frac12\frac{\omega^{1/2}+\omega^{-1/2}}{(\omega^{1/2} - \omega^{-1/2})^2}[\Sigma_n]^2,
\end{equation}
where $[\Sigma_n]^2$ denotes the self intersection number of $\Sigma_n$. 
When $p$ is odd, \eqref{eq:fn0-2} is valid with the sign if square roots are given by the rule \eqref{eq:sign}.

In order to calculate $k_j^l$ for other $l$, we note the following lemma.
\begin{Lemma}\label{lem:action}
Let $g\in G$, and the action of $g$ on the $\Spinc$-structure $c$ be denoted by $\hat{g}$. 
For a connection $A$ on $L$, if there exists $u\in \G_0$ which satisfies $\hat{g} A =uA$, i.e., $[A]\in J^G$,  then we can define another action $\hat{g}^\prime$ of $g$ on $c$ so that $\hat{g}^\prime A = A$.  
\end{Lemma}
\proof
Consider the action $(u^{-1}\circ\hat{g})$. 
Then  $(u^{-1}\circ\hat{g}) A = A$. 
In particular, we have $(u^{-1}\circ\hat{g})^p A = A$. 
Note that $(u^{-1}\circ\hat{g})^p$ is an element of $\G_0$. 
Therefore  $(u^{-1}\circ\hat{g})^p=1$
Thus $\hat{g}^\prime := (u)^{-1}\circ\hat{g}$ is a required action.
\endproof
Thus, for any $t_l=[A_l]\in J^G$, we can redefine the $G$-action on $c$ so that $A_l$ is $G$-invariant. 
Hence, $k_j^l$ are also calculated by the $G$-index formula. 
However, the contributions from fixed points for the redefined action are different from the original ones as 
\begin{equation}\label{eq:g-ind}
\ind_g D_{A_l} =\sum_{j=0}^{p-1} \zeta^j k_j^l = \sum_{n} \F^l_n(g), 
\end{equation}
where $\F_n^l(g)$ are calculated as in \eqref{eq:fn0} and \eqref{eq:fn0-2} for the redefined $g$ action on $c$.

For different $l_0$ and $l_1$, the difference between $\F_n^{l_0}(g)$ and $\F_n^{l_1}(g)$ is given as follows.
We can consider that a representation of $t_l\in J^G$ is given as  a triplet $(S^+_l,\phi_l, A_l)$ of a $G$-spinor bundle $S^+_l$, a trivialization $\phi_l$ at $x_0$, and a $G$-invariant connection $A_l$ on the determinant line bundle $L_l=\det S^+_l$. 
For $l_0$ and $l_1$, the difference between $(S^+_{l_0},\phi_{l_0}, A_{l_0})$ and $(S^+_{l_1}, \phi_{l_1}, A_{l_1})$ is given as a flat $G$-line bundle $\LL_{l_1 l_0}$:
$$
(S^+_{l_1}, \phi_{l_1}, A_{l_1}) = \LL_{l_1 l_0}\otimes(S^+_{l_0},\phi_{l_0}, A_{l_0}).
$$
For each component $X_n\subset X^G$, the weight of $g$-action on the fiber of $\LL_{l_1 l_0}$ at $x_n\in X_n$ is given as a complex number $\lambda_n^{l_1l_0}$, which is a $p$-th root of $1$. 
Then the relation between $\F_n^{l_0}(g)$ and $\F_n^{l_1}(g)$ is given as
\begin{equation}\label{eq:relation}
\F_n^{l_1}(g) = \lambda_n^{l_1l_0}\F_n^{l_0}(g).
\end{equation}

Before ending this subsection, we give a useful lemma for lifts of the $G$-action to a $\Spinc$-structure.
\begin{Lemma}\label{lem:lift}
Let $G=\Z_p$, and $X$ be a closed oriented $G$-manifold which has no $2$-torsion in $H_1(X;\Z)$. %the $1$st integral homology group, that is, if $2x=0$ for $x\in H_1(X;\Z)$, then $x=0$.
If the determinant line bundle of a $\Spinc$-structure $c$ on $X$ is $G$-equivariant, then the $G$-action lifts to $c$, that is, $c$ is $G$-equivariant or $G$-action is of even or odd type with respect to $c$ when $p=2$. 
\end{Lemma}
\proof
If there is no $2$-torsion in $H_1(X;\Z)$, then there is a bijective correspondence between the set of equivalence classes of $\Spinc$-structures and the set of equivalence classes of determinant line bundles. 
For $g\in G$, let $\bar{g}$ be the action of $g$ on $P_{\SO}\times_X P_{\U(1)}$. 
Consider the $2$-fold covering $P_{\Spinc}\to P_{\SO}\times_X P_{\U(1)}$. 
Since $\bar{g}^*P_{\Spinc}$ is isomorphic to $P_{\Spinc}$, we can lift $\bar{g}$ to $P_{\Spinc}$. 
Therefore the $G$-action on $P_{\SO}\times_X P_{\U(1)}$ lifts to a $\hat{G}$-action on $P_{\Spinc}$.  
\endproof
%
%
%%%%
\subsection{An example of application in the case when $G=\Z_2$.}
%%%% 
%
The next proposition which is an application of \thmref{thm:odd} is also a generalization of Fang's result. 
(Compare with Corollary $4$ of \cite{Fang}.)
However, this is not a ``new result'', for this can be proved by the adjunction inequality. (See \exref{ex:z2}.) 
Nevertheless, we state this as an example of application.
\begin{Proposition}\label{prop:adjunction}
Let $G=\Z_2$, and $X$ be a closed oriented $4$-dimensional $G$-manifold with $b_+\geq 2$ and $b_+^G\geq 1$, and suppose that $H_1(X;\Z)$ has no $2$-torsion. 
Suppose that there is a $\Spinc$-structure $c$ whose determinant line bundle is trivial, and $\SW_X (c) \not\equiv 0 \mod 2$. 
Let $d(c)$ be as in \eqref{eq:dim-moduli}.
If the $G$-action has no isolated fixed point,
then the following inequality holds:
\begin{align*}
1-b_1+b_+&\geq 2(1-b_1^G + b_+^G), \text{ when } d(c)\equiv 0 \mod 4,\\
1-b_1+b_+&\geq 2(-b_1^G + b_+^G), \text{ when } d(c)\equiv 2 \mod 4.
\end{align*} 
\end{Proposition}
\proof
Note that $c$ is the $\Spinc$-structure which is determined by a $\Spin$-structure. 
Since the determinant line bundle $L$ is trivial, we can define a $G$-action on $L$ which is the product of the $G$-action on $X$ and trivial action on fiber. 
Therefore the $G$-action lifts to $c$ by \lemref{lem:lift}. 
The lifted action may be of odd or even type with respect to $c$. 
We take the trivial flat connection $A_0$ on $L$ as the origin of the Jacobian torus $J$.
As is known widely, 
a $G$-action is of even type if and only if the fixed point set is isolated. 
On the other hand, a $G$-action is of odd type if and only if the fixed point set is $2$-dimensional.  
(See e.g. \cite{AB}.)
Therefore, if the $G$-action is of even type, then it must be free by the assumption.

Suppose that the $G$-action is of odd type. 
By the $G$-index formula (put $\omega=-1$ and $\nu_n=1$ in \eqref{eq:fn0-2}), we have $\F^0_n(g) =0$ for any component $X_n$ of $X^G$.
The relation \eqref{eq:relation} implies $\F^l_n(g) =0$ for any $l$ and $n$.
 
Therefore, we have $k^l_1 = k^l_3 = \frac12 \ind D_{A_0}$ for any $l$. 
By \thmref{thm:odd} with the assumption of mod $2$ non-vanishing of $\SW_X (c)$, it holds that, for any partition $(d_1, d_3)$ of $d(c)/2$, there exist $l$ and $j$ which satisfy
$$
2k_j^l \geq 2 d_j + 1 - b_1^G + b_+^G. 
$$  
Therefore we have
\begin{align*}
\ind D_{A_0} & \geq \frac{d(c)}2 + 1 - b_1^G + b_+^G, \hspace{1.33cm}\text{ when } d(c)\equiv 0 \mod 4,\\
\ind D_{A_0} & \geq \left(\frac{d(c)}2-1\right) + 1 - b_1^G + b_+^G, \text{ when } d(c)\equiv 2 \mod 4.
\end{align*}
On the other hand, from the formula of the dimension of the moduli \eqref{eq:dim-moduli}, we have
$$
\ind D_{A_0} = \frac12 (d(c) + 1 - b_1 + b_+).
$$
This formula with above two inequality implies the proposition.

In the even case, the $G$-action should be free. 
In the free case, the theorem is obvious from the Lefschetz formula and the $G$-signature formula.
\endproof
\begin{Example}\label{ex:z2}
Concrete examples of $G=\Z_2$-actions are given as follows.
Let $X$ be the $K3$ surface of Fermat type,
$
X=\{[z_0,z_1,z_2,z_3]\in \CP^3| z_0^4+z_1^4+z_2^4+z_3^4=0\}.
$
Let $G$ act on $X$ by the permutation of two coordinates. 
Then the fixed point set is a complex curve $C$ whose genus is $3$ and self-intersection number is $4$. 

Another example of $b_1>0$ is $4$-torus. 
Let $X$ be the direct product of two copies of $2$-torus.
Let $G$ act on the first $2$-torus by multiplication by $-1$, and on the second trivially. The fixed point set consists of four $2$-tori whose self-intersection number are $0$. 

Let us verify \propref{prop:adjunction} for these examples.
It is well-known that $\SW_X(c_0)=\pm 1$ for the $K3$ surface and the $4$-torus \cite{Taubes}. 
Note that, for a (V-)manifold $Y$, it holds that 
$$
1 - b_1(Y) + b_+(Y) = \frac12 (\chi (Y) + \Sign (Y)).  
$$
Therefore, by using the Lefschetz formula and the $G$-signature theorem, we have
\begin{align*}
1 - b_1^G + b_+^G &= \frac12 (\chi (X/G) + \Sign (X/G))\\
&= \frac12\left\{\frac12 (\chi(X) + \chi (C)) + \frac12(\Sign(X) + [C]^2)\right\}\\
&= \frac12\left\{\frac12 (\chi (X) + \Sign (X))\right\}\\
&= \frac12 (1 - b_1 + b_+).
\end{align*}
We use the adjunction formula at the third equality.
From this calculation, we see that the adjunction inequality $ \chi(C) +  [C]^2 \leq 0$ proves \propref{prop:adjunction}.
\end{Example}
\begin{Remark}
We can construct similar $G$-actions on homology $4$-tori obtained by the 'knot surgery' construction according to \cite{RS} and \cite{FS}.
(See also \exref{ex:z3}.)
\end{Remark}
%
%
%%%%
%
%
%%%%
\subsection{Examples of the case when $G=\Z_3$}\label{subsec:ex-z3}
%%%%
%
%
This subsection treats with the case when $G=\Z_3$. 
In the following, we assume that the $G$-action is {\it pseudofree}, that is, the $G$-action has only isolated fixed points.
In such a case, fixed points are classified into two types of representations:
\begin{itemize}
\item The type $(+)$: $(1,2)=(2,1)$.
\item The type $(-)$: $(1,1)=(2,2)$.
\end{itemize}
Let $m_+$ be the number of fixed points of the type $(+)$, and $m_-$ be that of the type $(-)$.

We give examples of pseudofree $G$-actions which imply the mod-3 vanishing of Seiberg-Witten invariants. 
\begin{Example}\label{ex:z3}
Let $X$ be the direct product of a $2$-torus and a Riemann surface of genus $3h$ ($h\geq 1$). 
We construct a $G$-action on $X$ as follows. 
Let us consider the lattice $\Z\oplus\zeta\Z\subset \C$, where $\zeta=\exp (2\pi\sqrt{-1}/3)$, and let $T_1$ be the  $2$-torus $\C_1/(\Z\oplus\zeta\Z)$ with a $G$-action, where the $G$-action is given by the multiplication by $\zeta$. 
%For each $j=0,1,2$, let us consider a $2$-torus $T_j=\C/(\Z\oplus\zeta\Z)$ with the $G$-action which is defined by the multiplication by $\zeta^j$. 
Next consider a $2$-sphere, and let $G$ act on the $2$-sphere by $2\pi/3$-rotation.
Taking a free point $q$ on the $2$-sphere, and glueing $3$ copies of a Riemann surface of genus $h$ to the $2$-sphere at three points $q$, $gq$, $g^2q$, we obtain a Riemann surface $\Sigma_{3h}$ of genus $3h$ with a $G$-action. 
Let $X$ be $T_1\times \Sigma_{3h}$ with the diagonal $G$-action. 

Now let us examine \thmref{thm:main}. 
First note that the fixed point set of $T_1$ consists of three points $p_0$, $p_1$ and $p_2$, and all of them have same type of representation: $T(T_1)_{p_n}\cong \C_1$. 
On the other hand, $\Sigma_{3h}$ have two fixed points $q_+$ and $q_-$, and they have opposite representations each other. 
(We assume that $q_+$ is the fixed point such that $T(\Sigma_{3h})_{q_+}\cong \C_2$. )
Therefore, $X$ has six fixed points, and three of them are of the type $(+)$, and the other three are of the type $(-)$.

Note that $\chi(X)=\Sign (X)=0$ and $X$ is spin. 
We take the $\Spinc$-structure $c_0$ which is determined by a $\Spin$-structure.
Note that $d(c_0)=0$. 
We consider the $G$-action on $c_0$ which induces the $G$-action on the determinant line bundle $L$ which is the product of the $G$-action on $X$ and the trivial action on fiber.
Take the trivial flat connection $A_0$ on $L$ as the origin of the Jacobian torus $J_X$.

The Jacobian $J_X$ is of the form $J_X=J_{T_1}\times J_{\Sigma_{3h}}$. 
For a fixed point $t=(a,b)\in J^G_X$, the corresponding flat $G$-bundle $\LL_t$ is written as $\LL_t = \pi_1^*\LL_a\otimes\pi_2^*\LL_b$, where $\pi_1$ (resp. $\pi_2$) is the projection to $T_1$ (resp. $\Sigma_{3h}$), and $\LL_a$ is the flat $G$-bundle on $T_1$ associated to $a\in J_{T_1}^G$ and $\LL_b$ is similar.

Now let us attempt to classify flat $G$-bundles on a Riemann surface. 
Temporarily, we consider more general situation that $G_p=\Z_p$ acts pseudofreely on a Riemann surface $\Sigma_g$ of genus $g$. 
Let $\{p_n\}$ be the fixed point set. 
Consider a divisor $D$ on $\Sigma_g$: $D=\sum_n d_n p_n$. 
Then we can construct a $G_p$-line bundle $\LL_D$ on $\Sigma_g$ which satisfy $\LL_D|_{p_n}\cong (T\Sigma_g|_{p_n})^{\otimes d_n}$. 
Note that $c_1(\LL_D)=\sum d_n$. 
In this situation, we can prove the following.
\begin{Proposition}\label{prop:flat}
Let $\LL$ be a $G_p$-line bundle on $\Sigma$ which satisfy $\LL|_{p_n}\cong (T\Sigma_g|_{p_n})^{\otimes d_n}$. 
Then $c_1(\LL)\equiv c_1(\LL_D) \mod p$.
\end{Proposition}
\proof
Let us consider the line bundle $\LL\otimes\LL_D^{-1}$. 
Then there is a line bundle $\bar{\LL}$ on $\Sigma_g/G_p$ which satisfies $\pi^* \bar{\LL}\cong \LL\otimes\LL_D^{-1}$, where $\pi\colon\Sigma_g\to\Sigma_g/G_p$ is the quotient map.
Noting that $c_1(\LL\otimes\LL_D^{-1})=\pi^*c_1(\bar{\LL})$, and $\pi^*\colon H^2(\Sigma_g/G_p;\Z)\to H^2(\Sigma_g;\Z)$ is multiplication by $p$, we have the proposition.
\endproof
Let us apply \propref{prop:flat} to $\Sigma_{3h}$ with the $G$-action. 
Since the fixed point set is $\{q_+,q_-\}$, the divisor $D$ is of the form $D = d_+q_+ + d_-q_-$. 
Since $\LL_b$ is trivial, we have $0=c_1(\LL_b) \equiv d_++d_-\mod 3$. 
Therefore, the following holds.
\begin{Lemma}\label{lem:L_b}
For any $b\in J_{\Sigma_{3h}}^G$, $\LL_b$ is isomorphic to $\LL_D$ such that $D= 0$ or $q_+-q_-$ or $2q_+-2q_-$.
\end{Lemma}
For $b\in J_{\Sigma_{3h}}^G$, let us denote the weight of the $G$-action on the fiber of $\LL_b$ at $q_+$ (resp. $q_-$) by $\lambda^b_+$ (resp. $\lambda^b_-$).
Similarly, for $a\in J_{T_1}^G$, denote the weight of $\LL_a$ at $p_i\in T_1^G$ by  $\lambda^a_i$.
Note that  $\F_{(p_i,q_\pm)}^{(0,0)}(g)$ for the origin $(0,0)\in J_X^G$ at $(p_i,q_\pm)\in X^G$ is given by $\F_{(p_i,q_\pm)}^{(0,0)}(g) = \pm\frac13$.
(See \eqref{eq:fn0}.) 
Therefore  $\F_{(p_i,q_\pm)}^{(a,b)}(g)$ for $(a,b)\in J_X^G$ at $(p_i,q_\pm)$ is written as 
$$
\F_{(p_i,q_\pm)}^{(a,b)}(g) = \pm\frac13\lambda_i^a\lambda_{\pm}^b.
$$
By \lemref{lem:L_b}, we have $\lambda^b_+=\lambda^b_-$.%= \zeta^k$, where $k$ is one of $0,1,2$.
Hence we obtain 
\begin{equation}\label{eq:fab}
\sum_{x\in X^G}\F^{(a,b)}_x(g) = \frac13\left(\sum_{i=0}^2\lambda_i^a\right)(\lambda_+^b-\lambda_-^b) = 0.
\end{equation}
Similarly we obtain 
\begin{equation}\label{eq:fab2}
\sum_{x\in X^G}\F^{(a,b)}_x(g^2) = 0,
\end{equation} 
for any $(a,b)\in J_X^G$.

By \eqref{eq:fab} and \eqref{eq:fab2}, the $G$-index formula for the Dirac operator of $t_l=[A_l]\in J^G$ is given as
\begin{align*}
\ind_g D_{A_l} &= k_0^l + \zeta k_1^l + \zeta^2 k_2^l = 0,\\
\ind_{g^2} D_{A_l} &= k_0^l + \zeta^2 k_1^l + \zeta k_2^l = 0,\\
\ind_1 D_{A_l} &= k_0^l + k_1^l + k_2^l = - \frac18 \Sign (X)=0.
\end{align*}
%where $\zeta=\exp (2\pi\sqrt{-1}/3)$.
Solving these equations, we obtain
$$
k_0^l =k_1^l=k_2^l=0.
$$

Now let us check that inequalities \eqref{eq:ineq} are satisfied. 
First let us compute $1-b_1^G+b_+^G$.
%From the Lefschetz formula and the $G$-signature formula, we have
The Lefschetz formula implies that
\begin{equation}\label{eq:G-euler}
\chi (X/G) = \frac13(\chi (X) + 2(m_+ +m_-)).
\end{equation}
On the other hand, the $G$-signature theorem (Cf.\cite{AB}) implies that
\begin{align}
\Sign (g,X) &= \Sign (g^2,X) = \frac13 (m_+ - m_-),\\\label{eq:G-sign}
\Sign (X/G) &=  \frac13 \left\{ \Sign(X) + \frac23 (m_+ - m_-)\right\}.
\end{align}
Since $\chi (X)=\Sign(X)=0$, we have,
\begin{equation}\label{eq:b3}
1-b_1^G+b_+^G = \frac12 (\chi(X/G)+\Sign(X/G))=\frac19(4m_+ + 2m_-)=2.
\end{equation}

Since the dimension of the moduli $d(c_0)$ is $0$, all $d_j$ in \eqref{eq:ineq} should be $0$. 
Therefore inequalities \eqref{eq:ineq} are satisfied as,
$$
2k_j^l = 0 < 2=1-b_1^G+b_+^G,
$$
for any $j,l$, and hence \thmref{thm:main} implies that  $\SW_X(c_0)\equiv 0\mod 3$.

On the other hand, we can calculate the Seiberg-Witten invariants of $X_g=T^2\times\Sigma_g$. 
The answer is given as follows: for the $\Spinc$-structure $c_0$ which is determined by a $\Spin$-structure, 
\begin{equation}\label{eq:SW_Xg}
\SW_{X_g}(c_0) = \pm
\begin{pmatrix}
2g-2\\
g-1
\end{pmatrix}.
\end{equation}
It is easy to see that this is divisible by $3$ if $g=3h$.
Thus, \thmref{thm:main} holds.

There are several methods to prove \eqref{eq:SW_Xg}.
One method is  Witten's calculation \cite[pp.786--792]{Witten}. 
The canonical divisor of $X_g$ is written as $c_1(K) = (2g-2)P.D.[T\times pt]$.
For a generic choice of $\eta\in H^0(X_g,K)$, a Seiberg-Witten solution corresponds to a factorization $\eta = \alpha \beta$, where $\alpha$ and $\beta$ are holomorphic sections of $K^{1/2}\otimes L^{\pm 1}$. 
Since $L$ of our case is trivial, the number of possibilities of factorizations $\eta = \alpha \beta$ coincides with the right hand side of \eqref{eq:SW_Xg}. 
Furthermore, we can see that all solutions have same sign also by \cite{Witten}.

An alternative way to prove \eqref{eq:SW_Xg} is as follows.
First consider $X_g$ as $S^1\times M$, where $M=S^1\times \Sigma_g$. 
Next determine the Seiberg-Witten invariants of $M$ by, for instance, Turaev torsion of $M$. 
Then use the formula $\SW_{S^1\times M}(\tilde{c})=\SW_{M}(c)$ where $\tilde{c}$ is the pull-back of $c$. 
When $g\geq 2$, Turaev torsion of $S^1\times \Sigma_g$ is written as $\pm (t-1)^{2g-2}$, where $t$ is the homology class represented by $S^1$, and $c_0$ corresponds to the term of order $g-1$. (See \cite[pp.93--96]{Turaev}.)
\end{Example}
\begin{Remark}\label{rem:z3}
Similar examples can be constructed via the 'knot surgery' construction of Fintushel and Stern \cite{FS}. 
Remove three copies of $T^2\times D^2$ from $X=T^2\times \Sigma_{3h}$ which are mapped to each other by the $G$-action, and denote the resulting manifold by $X^\prime$. 
According to \cite{FS}, let $K$ be a knot in $S^3$, and $E_K$ be the exterior. 
Then glueing $S^1\times E_K$ to each boundary of $X^\prime$ gives an example.  
This manipulation changes the Seiberg-Witten invariant by a multiple of $3$ \cite{FS}. 
\end{Remark}
\begin{Remark}
We can construct an example of $G$-action such that the Seiberg-Witten invariant {\it does not} vanish modulo $3$ and there exists $l$ for which \eqref{eq:ineq} {\it does not} hold. 
Let $T_i$ be the  $2$-torus $\C_i/(\Z\oplus\zeta\Z)$ with the $G$-action given by the multiplication by $\zeta^i$ ($i=1,2$). 
Remove a small $G$-invariant neighborhood of a fixed point of each $T_i$. 
Since fixed points of $T_1$ and $T_2$ have opposite representations, we can glue their boundaries $G$-equivariantly, and the resulting manifold is a Riemann surface $\Sigma_2$ of genus $2$ with a $G$-action whose fixed point set consists of four points. 
Now consider the $4$-manifold $T_1\times \Sigma_2$ with the diagonal $G$-action. 
Then we can prove that this is a required example.
\end{Remark}

\end{document}